\documentclass[12pt,reqno,openbib,runningheads,a4paper]{amsart}%
\usepackage{graphicx}
\usepackage{amssymb}
\usepackage{amsfonts}
\usepackage{amsmath}
\usepackage{graphicx}
\usepackage{lineno}
\usepackage[authoryear]{natbib}
\usepackage[dvipsnames]{xcolor}
\usepackage{epsf}
\usepackage{lscape}
\usepackage[legalpaper,bookmarks=blue,filecolor=blue,linkcolor=blue,citecolor=blue,colorlinks=true]%
{hyperref}%
\providecommand{\U}[1]{\protect \rule{.1in}{.1in}}
\providecommand{\U}[1]{\protect \rule{.1in}{.1in}}
\makeatletter
\renewcommand{\@biblabel}[1]{}
\makeatother

\begin{document}

\begin{center}
{\Large \textbf{Extreme value theory based confidence intervals for the
parameters of a symmetric L\'{e}vy-stable distribution}}\medskip \medskip

{\large Djamel Meraghni}$^{\ast}${\large , Louiza Soltane}\medskip

{\small \textit{Laboratory of Applied Mathematics, Mohamed Khider University,
Biskra, Algeria}}\bigskip \medskip
\end{center}

\noindent \textbf{Abstract}\medskip

\noindent We exploit the asymptotic normality of the extreme value theory
(EVT) based estimators of the parameters of a symmetric{\Large \textbf{ }%
}L\'{e}vy-stable distribution, to construct confidence intervals. The accuracy
of these intervals is evaluated through a simulation study.\medskip

\medskip


\vfill

\noindent{\small $^{\text{*}}$Corresponding author\newline \noindent
\textit{E-mail addresses :}\newline \texttt{djmeraghni@yahoo.com}
(D.~Meraghni)\newline \texttt{louiza\_stat@yahoo.com} (L.~Soltane)}

\section{\textbf{Introduction\label{sec1}}}

\subsection{L\'{e}vy-stable Distributions\-}

\noindent

\noindent The L\'{e}vy-stable distribution, also called stable, $\alpha
$-stable or stable Paretian, represents a rich class of probability
distributions. Introduced in 1920's by \cite{Levy25}, while investigating the
behavior of normalized sums of independent identically distributed (iid)
random variables (rv's), it has got an increased attention in the last decades
for at least two good reasons. First, it is theoretically supported by the
generalized central limit theorem which states that the $\alpha$-stable law is
the only possible limit distribution for properly normalized and centered sum
of iid rv's. Second, it allows skewness and fat tails meaning that it is
suitable for data collected in areas as diverse as finance, hydrology,
meteorology,... Indeed, a great deal of empirical evidence indicates that
these data can be so heavy-tailed that they are poorly described by the
largely used Gaussian distribution. In other words, the stable model provides
a much better fit for heavy-tailed observations sets than the commonly adopted
normal one does.\smallskip

\noindent The extreme value theory (EVT), which proved to be an excellent tool
in risk management, could be applied to estimate the parameters characterizing
a stable distribution in order to determine the appropriate model for a given
data set. In the sequel, let $\overset{d}{=},$ ${\normalsize \overset
{p}{\rightarrow}}\ $and $\overset{d}{\rightarrow}$ stand for equality in
distribution, convergence in probability and convergence in distribution
respectively and let $\mathcal{N}\left(  m,v^{2}\right)  $ denote the normal
distribution with mean $m\in \mathbb{R}$ and variance $v^{2}>0.$\smallskip

\noindent A rv $X$\ is said to be L\'{e}vy-stable if and only if, for
$n\geq2,$\ $\exists$ $a_{n}>0,$\ $b_{n}\in \mathbb{R}$\ such that%
\[
\frac{\left(  X_{1}+...+X_{n}\right)  -b_{n}}{a_{n}}\overset{d}{=}X,
\]
where $X_{1},...,X_{n}$\ are independent copies of $X.$ It is shown that
$\exists$ $0<\alpha \leq2$ such that $a_{n}=n^{1/\alpha},$ (see, e.g.,
\citeauthor{Feller71}, \citeyear{Feller71}).\smallskip

\noindent Except from three special cases, a stable rv suffers from the lack
of closed-form expressions for its distribution function (df) and probability
density function (pdf). However, it is typically described by its
characteristic function $\varphi$ which has many representations. The most
famous one is defined for $t\in \mathbb{R}$ by%

\[
\varphi \left(  t\right)  =\left \{
\begin{array}
[c]{lcc}%
\exp \left \{  i\mu t-\sigma^{\alpha}\left \vert t\right \vert ^{\alpha}\left(
1-i\beta sign\left(  t\right)  \tan \frac{\alpha \pi}{2}\right)  \right \}  &
\text{\textit{for}} & \alpha \neq1,\\
\exp \left \{  i\mu t-\sigma \left \vert t\right \vert \left(  1+i\beta sign\left(
t\right)  \frac{2}{\pi}\log \left \vert t\right \vert \right)  \right \}  &
\text{\textit{for}} & \alpha=1,
\end{array}
\right.
\]
where%
\[
i^{2}=-1\text{ and }sign(t):=\left \{
\begin{array}
[c]{ccc}%
1 & \text{if} & t>0,\\
0 & \text{if} & t=0,\\
-1 & \text{if} & t<0.
\end{array}
\right.
\]
As we may see, this family of distributions is characterized by four
parameters :

\begin{itemize}
\item[$\bullet$] $0<\alpha \leq2:$ stability index, tail exponent or shape
parameter.\smallskip

\item[$\bullet$] $\sigma>0:$ scale parameter.\smallskip

\item[$\bullet$] $-1\leq \beta \leq1:$ skewness parameter.\smallskip

\item[$\bullet$] $\mu \in \mathbb{R}:$ location parameter.
\end{itemize}

\noindent Using a notation of \cite{ST94}, a rv $X$ with stable distribution
will be written as $X\sim S_{\alpha}(\sigma,\beta,\mu).$ The three cases where
we have explicit formulas for the pdf are the very popular Gaussian
distribution $S_{2}(\sigma,0,\mu)$ and the lesser known models of Cauchy
$S_{1}(\sigma,0,\mu)$ and L\'{e}vy $S_{1/2}(\sigma,1,\mu).$ The tail exponent
$\alpha,$ which is the most important among all four parameters, indicates the
rate at which the tails of the distribution taper off. For $0<\alpha<2,$ the
$k$th $\left(  k=1,2,...\right)  $ moment of a stable rv is finite if and only
if $k<\alpha,$ whereas for $\alpha=2$ all the moments exist. In particular,
the distribution mean only exists when $1<\alpha \leq2$ and is equal to the
location parameter $\mu.$ For $0<\alpha<2,$ the variance is infinite and the
distribution tails are asymptotically equivalent to those of a
Pareto\ distribution, i.e., they exhibit a power-law behavior.

\subsection{Heavy Tails Property of $S_{\alpha}(\sigma,\beta,\mu)$}

\begin{figure}
[ptb]
\begin{center}
\includegraphics[
height=2.7838in,
width=5.1794in
]%
{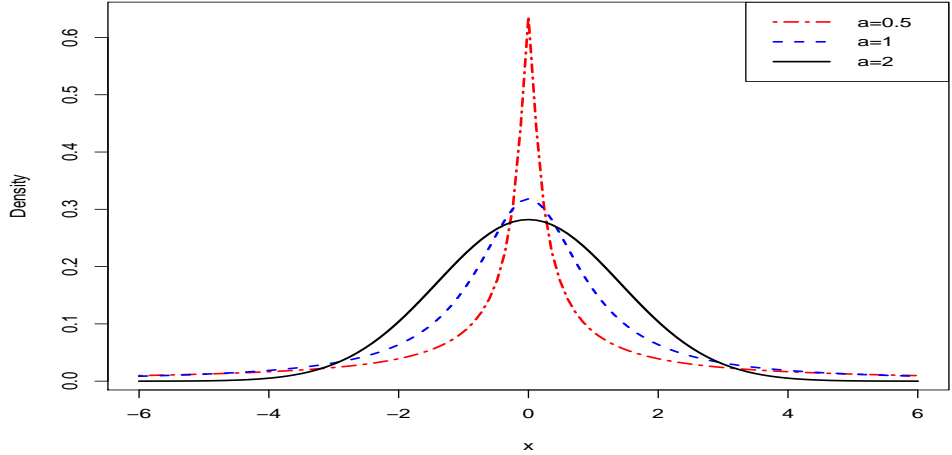}%
\caption{Stable densities for different values of $\alpha$ with $\beta=\mu=0$
and $\sigma=1.$ The solid line corresponds to the Gaussian model.}%
\label{F1}%
\end{center}
\end{figure}

\noindent In general, the upper and lower tails of a L\'{e}vy-stable
distribution asymptotically exhibit a Pareto-like behavior, i.e. they fall off
like a power function. The rate of decay is governed by the stability index :
the smaller $\alpha,$ the slower the decay and hence the heavier the
distribution tails, as shown in \autoref{F1}.%

More precisely, for a rv $X\sim S_{\alpha}(\sigma,\beta,\mu),$ the following
result holds (see e.g., \citeauthor{ST94}, \citeyear{ST94}, page $16$).%

\begin{equation}
\underset{x\rightarrow \infty}{\lim}x^{\alpha}P\left(  X>x\right)  =C_{\alpha
}\frac{1+\beta}{2}\sigma^{\alpha}\text{ and }\underset{x\rightarrow \infty
}{\lim}x^{\alpha}P\left(  X<-x\right)  =C_{\alpha}\frac{1-\beta}{2}%
\sigma^{\alpha}, \label{P}%
\end{equation}
where%
\begin{equation}
C_{\alpha}:=\left(  \int \nolimits_{0}^{\infty}x^{-\alpha}\sin xdx\right)
^{-1}=\frac{2}{\pi}\Gamma \left(  \alpha \right)  \sin \frac{\pi \alpha}{2},
\label{C}%
\end{equation}
with $\Gamma$ being the gamma function defined, for $u>0,$ by $\Gamma \left(
u\right)  =\int \nolimits_{0}^{\infty}x^{u-1}e^{-x}dx.$ From equations $\left(
\ref{P}\right)  ,$ we get what is specifically called tail balance conditions.
That is, we have, as $x\rightarrow \infty,$%

\begin{equation}
\frac{P\left(  X>x\right)  }{P\left(  \left \vert X\right \vert >x\right)
}\rightarrow \frac{1+\beta}{2}=:p\text{ and }\frac{P\left(  X<-x\right)
}{P\left(  \left \vert X\right \vert >x\right)  }\rightarrow \frac{1-\beta}%
{2}=:q=1-p, \label{balanceconditions}%
\end{equation}
Let $F$ and $G$ denote the df's of $X$ $\sim$ $S_{\alpha}(\sigma,\beta,\mu)$
and $Z=\left \vert X\right \vert $ respectively. It is obvious that $F$ and $G$
are related by%
\[
G(x)=F(x)-F(-x),\text{ }x>0.
\]
From relation $\left(  \ref{P}\right)  ,$ we get that the distribution tail of
$Z$ satisfies%
\begin{equation}
1-G(x)\sim C_{\alpha}\sigma^{\alpha}x^{-\alpha},\text{ as }x\rightarrow \infty,
\label{tail}%
\end{equation}
and%
\begin{equation}
\underset{t\rightarrow \infty}{\lim}\frac{1-G(tx)}{1-G(t)}=x^{-\alpha},\text{
}x>0. \label{RV1}%
\end{equation}
The latter means that $1-G$ is regularly varying at infinity with \ index
$-\alpha<0.$ For full details on regular variation, see, for instance,
Appendix B in \cite{deHF06}. From \cite{Gnedenko43}\textbf{,} relation
$\left(  \ref{RV1}\right)  $ is equivalent to say that $G$ is in Fr\'{e}chet
maximum domain of attraction. More precisely, for a sample $Z_{1},...,Z_{n}$
$(n\geq1)$ from the rv $Z,$ we have%
\[
\frac{\max \left(  Z_{1},...,Z_{n}\right)  }{G^{-1}\left(  1-n^{-1}\right)
}\overset{d}{\rightarrow}\Phi_{\alpha},
\]
$G^{-1}\left(  u\right)  :=\inf \left \{  x\in \mathbb{R},\text{ }G(x)\geq
u\right \}  ,$ $0<u<1,$ is the generalized inverse or quantile function of $Z$
and%
\[
\Phi_{\alpha}(x):=\left \{
\begin{array}
[c]{ll}%
\exp(-x^{-\alpha}), & x>0\\
0, & x\leq0
\end{array}
.\right.
\]
For further details and a complete description of this class of distributions,
we refer to the textbooks of \cite{Feller71}, \cite{Zolotarev86}, \cite{ST94}
and \cite{Nolan01}. On the other hand, there are available some very useful
computer programs, such that "STABLE", "Xplore" and the package "stabledist"
of the statistical software R (\citeauthor{IG96}, \citeyear{IG96}), specially
developed for numerical purposes (computing stable df's and pdf's, generating
stable rv's, estimating stable parameters,...).\smallskip

\noindent In this work, we concentrate on the case where $\beta=0,$ that is
when the distribution is symmetric about $\mu.$ In this case, the
characteristic function and the tail balance conditions respectively reduce to
the simpler forms%
\[
\varphi \left(  t\right)  =\exp \left \{  -\sigma^{\alpha}\left \vert t\right \vert
^{\alpha}+i\mu t\right \}  ,\text{ }t\in \mathbb{R},
\]
and%

\[
\underset{x\rightarrow \infty}{\lim}\frac{P\left(  X>x\right)  }{P\left(
\left \vert X\right \vert >x\right)  }=\underset{x\rightarrow \infty}{\lim}%
\frac{P\left(  X<-x\right)  }{P\left(  \left \vert X\right \vert >x\right)
}=\frac{1}{2}.
\]
\smallskip

\noindent The rest of the paper is organized as follows. \autoref{sec2}, is
devoted to a brief reminder on EVT-based estimators of the stable parameters.
In \autoref{sec3}, we use the asymptotic normality property of the estimators
to build confidence intervals for parameters $\alpha,$ $\mu$ and $\sigma.$
Finally, the accuracy of such intervals is investigated in a simulation study
in \autoref{sec4}.

\section{EVT-based estimation\textbf{\label{sec2}}}

\noindent The lack of explicit forms for the df and pdf severely hampers the
estimation of the distribution parameters. Nevertheless, several numerical
procedures of estimation based on the sample quantiles, the sample
characteristic function and maximum likelihood\ approaches, are proposed in
the literature. In a comparative study \cite{Ojeda01} notices that maximum
likelihood based methods are the most accurate but the slowest of all others.
On the other hand, the nature of the L\'{e}vy-stable distribution tails
suggests that EVT could play a major role in estimating its parameters. EVT is
a classical topic in probability theory and mathematical statistics, developed
for the estimation of occurrence probability of rare events. It permits to
extrapolate the behavior of distribution tails from the largest observed data.
EVT techniques have proven to be very useful where estimation of tail-related
quantities such as extreme value index, high quantiles, small exceedance
probabilities and mean excess function, is needed. The domains of application
of EVT\ include insurance (premium computation, large losses,...), finance
(asset returns, exchange rate,...), hydrology (floods, drought,...),
meteorology (extreme weather conditions,...), ecology (pollution peaks,...),
telecommunications (network traffic,...), physics (nuclear reactions,...).
EVT-based estimation approach has at least three advantages. It focuses only
on tail behavior and does not assume a parametric form for the entire
distribution. It provides estimators of explicit forms making estimate
computation easier and more direct. Finally, it produces estimators which
enjoy the asymptotic normality property leading to the construction of
confidence bounds for the unknown parameters. A very good variety of textbooks
may be consulted for a review of this topic and its multiple applications. We
can cite, for instance, \cite{deHF06}, \cite{EKM97}, \cite{RT97} and
\cite{BeGeS04}.\smallskip

\subsection{Estimating the Stability Index}

\noindent

\noindent The characteristic exponent $\alpha$ is the main parameter as it
governs the behavior of the distribution tails. Many estimators are proposed
for $\alpha$ via the EVT approach, among which the most popular is that
introduced by Hill (\citeauthor{Hill75}, \citeyear{Hill75}) as follows :%

\begin{equation}
\widehat{\alpha}_{n}=\widehat{\alpha}_{n}(k):=\left(  \frac{1}{k}%
{\displaystyle \sum_{i=1}^{k}}
\log Z_{n-i+1:n}-\log Z_{n-k:n}\right)  ^{-1}, \label{Hill}%
\end{equation}
where $Z_{1:n}\leq...\leq Z_{n:n}$ are the order statistics pertaining to a
sample $\left(  Z_{1},...,Z_{n}\right)  ,$ $n\geq1,$ from the rv $Z$ and
$k=k(n)$ is an integer sequence such that%
\begin{equation}
k\rightarrow \infty \text{ and }k/n\rightarrow0\text{ as }n\rightarrow \infty.
\label{k}%
\end{equation}
The consistency of $\widehat{\alpha}_{n}$ is proved in \cite{Mas82}, while its
almost sure convergence is established in \cite{Necir06a}. For the asymptotic
normality of $\widehat{\alpha}_{n}$ (and other related estimators), it is
required an additional assumption, known as the second-order condition of
regular variation (see \citeauthor{deHS96}, \citeyear{deHS96}), which
specifies the rate of convergence in $\left(  \ref{RV1}\right)  .$ That is, we
assume that there exist a constant $\rho<0,$ called second-order parameter,
and a function $A$ tending to zero and not changing sign near infinity, such
that for any $x>0,$ we have%
\begin{equation}
\underset{t\rightarrow \infty}{\lim}\frac{\left(  1-G(tx)\right)  /\left(
1-G(t)\right)  -x^{-\alpha}}{A(t)}=x^{-\alpha}\frac{x^{\alpha \rho}-1}%
{\rho/\alpha}. \label{RV2}%
\end{equation}
Note that when $1<\alpha<2,$ the condition $\left(  \ref{RV2}\right)  $ is
fulfilled. Indeed, using the expansion (to the second order) given in top of
page 95 in \cite{Zolotarev86}, yields that $G$ belongs to Hall's class of
heavy-tailed distributions (\citeauthor{Hall82}, \citeyear{Hall82}), which in
turn implies that $\left(  \ref{RV2}\right)  $ holds. A df $K$ is said to
belong Hall's class if%
\begin{equation}
1-K\left(  x\right)  =cx^{-1/\gamma}\left(  1+dx^{\rho/\gamma}+o\left(
x^{\rho/\gamma}\right)  \right)  ,\text{ as }x\rightarrow \infty, \label{Hall}%
\end{equation}
where $\gamma>0,\  \rho \leq0,$ $c>0,$ and $d\neq0.$ Hall's class, which is a
subset of the more general family of models with second-order regularly
varying tails, includes distributions (Burr, Fr\'{e}chet,...) that are most
commonly used in extreme event modelling. Among the works on the asymptotic
normality of $\widehat{\alpha}_{n},$ we can cite that of \cite{Peng98} who
proved that, if $\left(  \ref{RV2}\right)  $ holds, then for an integer
sequence $k$ satisfying $\left(  \ref{k}\right)  $ and $\lim_{n\rightarrow
\infty}\sqrt{k}A(n/k)=\lambda,$ with $\lambda$ finite, then
\begin{equation}
\sqrt{k}\left(  \widehat{\alpha}_{n}^{-1}-\alpha^{-1}\right)  \overset
{d}{\rightarrow}\mathcal{N}\left(  \lambda/(1-\rho),\alpha^{-2}\right)
,\text{ as }n\rightarrow \infty. \label{alpha-norm}%
\end{equation}
\cite{Weron01} discussed the performance of Hill's estimator $\widehat{\alpha
}_{n}$ and noted that for $\alpha \leq1.5$ the estimation is quite reasonable
but as $\alpha$ approaches $2,$ there is a significant overestimation when
considering samples of typical size (for an illustration, see \autoref{F2} and
\autoref{Tab1}). For such values of $\alpha,$ a very large number of
observations (a million or more) is needed in order to obtain acceptable
estimates and avoid misleading inference on the stability index, because the
true heavy tail nature of the distribution is visible only for extremely large
datasets. Fortunately, this kind of datasets are available nowadays and their
storage and treatment are made possible thanks to a very sophisticated
technology.\smallskip

\noindent The behavior of Hill's estimator (and therefore that of EVT-based
estimators) is affected by the number $k$ of upper order statistics to be used
in estimate computations. One needs to locate where the distribution tails
really begin because using too many data results in a big bias and too few
observations lead to a substantial variance. Consequently, one has to make a
trade-off between bias and variance in order to get an accurate estimate. To
this end, it seems reasonable that minimizing the mean squared error allows
for a compromise between the bias and variance components. On the other hand,
there exist several algorithms and data-adaptive procedures for the selection
of the optimal sample fraction of extreme values that guarantees the best
possible estimate (see, for instance, \citeauthor{CP01}, \citeyear{CP01},
\citeauthor{DHP01}, \citeyear{DHP01}, \citeauthor{FV04}, \citeyear{FV04} and
\citeauthor{NFA04}, \citeyear{NFA04}).

\subsection{Estimating the Location Parameter}

\noindent The empirical mean $\overline{X}:=n^{-1}\sum_{i=1}^{n}X_{i},$ which
is the natural estimator of the mean, is, in virtue of the central limit
theorem, asymptotically normal provided that the second moment is finite.
However, for $X\sim S_{\alpha}(\sigma,\beta,\mu)$ with $1<\alpha<2,$ the
latter theorem is not applicable because the variance of $X$ is infinite.
Therefore, the asymptotic normality of the sample mean $\overline{X}$ is not
established. To solve this problem, \cite{Peng01} proposed an asymptotically
normal estimator $\widehat{\mu}_{n}$ for $\mu,$ based on the the order
statistics $X_{1:n}\leq...\leq X_{n:n}$ associated to a sample $\left(
X_{1},...,X_{n}\right)  $ from $X,$ as follows :%
\[
\widehat{\mu}_{n}=\widehat{\mu}_{n}(k):=\widehat{\mu}_{n}^{(1)}+\widehat{\mu
}_{n}^{2}+\widehat{\mu}_{n}^{(3)},
\]
where%
\[
\widehat{\mu}_{n}^{2}=\widehat{\mu}_{n}^{2}(k):=\frac{1}{n}%
{\displaystyle \sum \limits_{i=k+1}^{n-k}}
X_{i:n}\text{ (trimmed mean),}%
\]%
\[
\widehat{\mu}_{n}^{(1)}=\widehat{\mu}_{n}^{(1)}(k):=\dfrac{k}{n}X_{k:n}%
\dfrac{\widehat{\alpha}_{n}^{(1)}}{\widehat{\alpha}_{n}^{(1)}-1}\text{ and
}\widehat{\mu}_{n}^{(3)}=\widehat{\mu}_{n}^{(3)}(k):=\dfrac{k}{n}%
X_{n-k+1:n}\dfrac{\widehat{\alpha}_{n}^{(3)}}{\widehat{\alpha}_{n}^{(3)}-1},
\]
with%
\[
\widehat{\alpha}_{n}^{(1)}=\widehat{\alpha}_{n}^{(1)}(k):=\left(  \frac{1}{k}%
{\displaystyle \sum \limits_{i=1}^{k}}
\log \left(  -X_{i:n}\right)  -\log \left(  -X_{k:n}\right)  \right)  ^{-1},
\]
and%
\[
\widehat{\alpha}_{n}^{(3)}=\widehat{\alpha}_{n}^{(3)}(k):=\left(  \frac{1}{k}%
{\displaystyle \sum \limits_{i=1}^{k}}
\log X_{n-i+1:n}-\log X_{n-k:n}\right)  ^{-1},
\]
being consistent estimators of $\alpha$ as well. The strong limiting behavior
of $\widehat{\mu}_{n}$ is studied in \cite{Necir06b} when constructing a
nonparametric sequential test with power $1$ for $\mu.$ For the asymptotic
normality of $\widehat{\mu}_{n},$ we notice that, by the expansion (to the
second order) and the relationship between the tails of $X$ respectively given
in pages 95 and 65 of \cite{Zolotarev86}, both tails of $F$ satisfy the
definition of Hall's model $\left(  \ref{Hall}\right)  .$ \cite{Peng01} proved
that, with a suitable choice of $k,$%
\begin{equation}
\frac{\sqrt{n}}{\tau(k/n)}\left(  \widehat{\mu}_{n}-\mu \right)  \overset
{d}{\rightarrow}\mathcal{N}(0,\delta^{2}),\text{ as}\;n\rightarrow \infty,
\label{mu-norm}%
\end{equation}
where%
\[
\delta^{2}:=1+\left(  \frac{\left(  2-\alpha \right)  \left(  2\alpha
^{2}-2\alpha+1\right)  }{2\left(  \alpha-1\right)  ^{4}}+\frac{\left(
2-\alpha \right)  }{\left(  \alpha-1\right)  }\right)  ,
\]
and%
\[
\tau^{2}(s):=\int_{s}^{1-s}\int_{s}^{1-s}\left(  u\wedge v-uv\right)
dF^{-1}(u)dF^{-1}(v),\text{ }0<s<1.
\]
It is shown in \cite{Peng01} that, as $n\rightarrow \infty,$%
\[
\sqrt{k/n}F^{-1}(k/n)\tau(k/n)\overset{P}{\rightarrow}-\left(  \frac{2-\alpha
}{2\left(  p^{2/\alpha}+(1-p)^{2/\alpha}\right)  }\right)  ^{1/2}%
(1-p)^{1/\alpha},
\]
which in our case $(\beta=0),$ may be rewritten into%
\begin{equation}
\tau(k/n)\sim-\frac{2\sqrt{k/n}F^{-1}(k/n)}{\sqrt{2-\alpha}},\text{ as
}n\rightarrow \infty. \label{sigmaapproxim}%
\end{equation}

\subsection{Estimating the Scale Parameter}

\noindent By combining relations $(\ref{tail})$ and $(\ref{C}),$ with some
approximations, \cite{MN07} provided a consistent estimator $\widehat{\sigma
}_{n}$ to the scale parameter $\sigma$ as follows :%
\[
\widehat{\sigma}_{n}:=Z_{n-k:n}\left(  \frac{k\pi}{2n\Gamma \left(
\widehat{\alpha}_{n}\right)  \sin \dfrac{\pi \widehat{\alpha}_{n}}{2}}\right)
^{1/\widehat{\alpha}_{n}},
\]
and proved that, with an adequate sequence $k,$%
\begin{equation}
\frac{\sqrt{k}}{\log(k/n)}\left(  \log \widehat{\sigma}_{n}-\log \sigma \right)
\overset{d}{\rightarrow}\mathcal{N}\left(  \lambda/(1-\rho),\alpha
^{-2}\right)  ,\text{ as }n\rightarrow \infty. \label{sigma-norm}%
\end{equation}

\section{Confidence bounds\textbf{\label{sec3}}}

\noindent Let us fix the confidence level of estimation to be $0<1-a<1$ and
let $z_{a}$ denote the $(1-a)$-quantile of the standard Gaussian distribution.
The first step, in the process of confidence interval construction, is to
determine the optimal sample fraction, that we denote by $k^{\ast},$ of
extreme observations involved in estimate computation. To this end, we adopt
the methodology of \cite{NFA04} who discussed and evaluated the performance of
the procedure proposed by \cite{RT97}. The latter consists in taking as
optimal the value of $k$ that minimizes%
\begin{equation}
RT(k):=\frac{1}{k}%
{\displaystyle \sum_{i=1}^{k}}
i^{\theta}\left \vert \widehat{\alpha}_{n}(i)-med(\widehat{\alpha}_{n}\left(
1\right)  ,...,\widehat{\alpha}_{n}\left(  k\right)  )\right \vert , \label{RT}%
\end{equation}
where $med$ stands for the median and $0\leq \theta \leq1/2.$ In other words, we
have%
\[
k^{\ast}:=\arg \min_{k}RT(k).
\]
Since we will be interested in the range $0<1/2<1/\alpha<1,$ then we choose
$\theta=0.3$ as indicated in \cite{NFA04}.\smallskip

\noindent The second step is to compute the estimate values which correspond
to the optimal number $k^{\ast}.$ Note that, for parameters $\alpha$ and
$\sigma,$ we only use the top observations in the $Z$-sample, whereas for
$\mu,$ we use the whole $X$-sample. Finally, we exploit the asymptotic
normality results $\left(  \ref{alpha-norm}\right)  ,$ $\left(  \ref{mu-norm}%
\right)  $ and $\left(  \ref{sigma-norm}\right)  $ to get asymptotic
confidence bounds for $\alpha,$ $\mu$ and $\sigma$ respectively. If we set
$\widehat{\alpha}^{\ast}:=\widehat{\alpha}_{n}(k^{\ast}),$ $\widehat{\mu
}^{\ast}:=\widehat{\mu}_{n}(k^{\ast})$ and $\widehat{\sigma}^{\ast}%
:=\widehat{\sigma}_{n}(k^{\ast}),$ then the respective $(1-a)\times
100\%$-confidence intervals for parameters $\alpha,$ $\mu$ and $\sigma$ are%
\[
\left(  \frac{\widehat{\alpha}^{\ast}}{1+z_{a/2}/\sqrt{k^{\ast}}}\text{
},\text{ }\frac{\widehat{\alpha}^{\ast}}{1-z_{a/2}/\sqrt{k^{\ast}}}\right)  ,
\]%
\[
\left(  \widehat{\mu}^{\ast}-z_{a/2}\frac{\widehat{\delta}^{\ast}\widehat
{\tau}^{\ast}}{\sqrt{n}}\text{ },\text{ }\widehat{\mu}^{\ast}+z_{a/2}%
\frac{\widehat{\delta}^{\ast}\widehat{\tau}^{\ast}}{\sqrt{n}}\right)  ,
\]
and%
\[
\left(  \exp \left \{  \log \widehat{\sigma}^{\ast}+z_{a/2}\frac{\log(k^{\ast
}/n)}{\widehat{\alpha}^{\ast}\sqrt{k^{\ast}}}\right \}  \text{ },\text{ }%
\exp \left \{  \log \widehat{\sigma}^{\ast}-z_{a/2}\frac{\log(k^{\ast}%
/n)}{\widehat{\alpha}^{\ast}\sqrt{k^{\ast}}}\right \}  \right)  ,
\]
where%
\[
\widehat{\delta}^{\ast}:=\left(  1+\frac{\left(  2-\widehat{\alpha}^{\ast
}\right)  \left(  2\widehat{\alpha}^{\ast2}-2\widehat{\alpha}^{\ast}+1\right)
}{2\left(  \widehat{\alpha}^{\ast}-1\right)  ^{4}}+\frac{\left(
2-\widehat{\alpha}^{\ast}\right)  }{\left(  \widehat{\alpha}^{\ast}-1\right)
}\right)  ^{1/2}\text{ and }\widehat{\tau}^{\ast}:=-\frac{2\sqrt{k^{\ast}%
/n}X_{k^{\ast}:n}}{\sqrt{2-\widehat{\alpha}^{\ast}}}.
\]

\section{Simulation study\textbf{\label{sec4}}}

\noindent We carry out a simulation study, by means of the statistical
software R (\citeauthor{IG96}, \citeyear{IG96}), to illustrate the finite
sample behaviors of the three estimators $\widehat{\alpha}_{n},$ $\widehat
{\mu}_{n}$ and $\widehat{\sigma}_{n}$ by computing their absolute biases (abs
bias) and mean squared errors (mse). We also evaluate the accuracy of the
confidence intervals (conf int) through their lengths and coverage
probabilities (cov prob). But first, we start by graphically checking Weron's
note (\citeauthor{Weron01}, \citeyear{Weron01}) on Hill's estimator of the
stability index. It is noteworthy that, for each experiment, we make $1000$
replications then we take our overall results by averaging over all the
individual results obtained at the end of each repetition.\smallskip

\begin{figure}
[ptb]
\begin{center}
\includegraphics[
height=2.8158in,
width=4.4572in
]%
{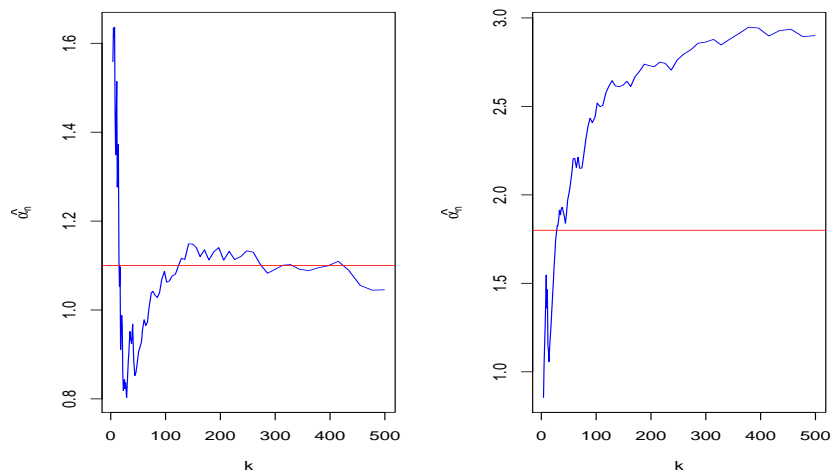}%
\caption{Plots of Hill's estimator (based on samples of size $5000$) of the
stability index $\alpha$ vs. the number $k$ of upper order statistics for
$\alpha=1.1$ (left) and $\alpha=1.8$ (right). The horizontal line represents
the true value of $\alpha.$}%
\label{F2}%
\end{center}
\end{figure}

\noindent We see on the right graph of \autoref{F2}, which is based on samples
of size $5000$ from stable distributions $S_{1.1}\left(  1,0,0\right)  $ and
$S_{1.8}\left(  1,0,0\right)  ,$ that there is no intermediate number $k$
which gives a good estimate for $\alpha$ and that estimates can even be above
the L\'{e}vy-stable regime. On the other side, the left panel shows that for
$\alpha=1.1,$ accurate estimates could be obtained for $k\in \{150,...,500\}.$%
\smallskip

\noindent For the estimation of the shape parameter $\alpha,$ we generate
$3000$ observations\ of symmetric $\alpha$-stable distributions $S_{\alpha
}\left(  \sigma,0,0\right)  $ with several values for parameters $\alpha$ and
$\sigma.$ The results are summarized in \autoref{Tab1}, where we note that, as
expected, the smaller the parameter values, the better the estimation. The
bottom of the table (corresponding to $\alpha=1.8$) shows that the estimation
is very poor for large $\alpha$-values and confirms the graphical conclusion
we made about the irrelevance of Hill's estimator, for large stability
indices, when built on the basis of datasets of typical sizes. For this
reason, only the values of $\alpha$ that are less than or equal to $1.5$ will
be considered thereafter. We gather the simulation results in \autoref{Tab2}
for the location parameter $\mu$ and in \autoref{Tab3} for the scale parameter
$\sigma.$ The former shows that the more $\alpha$ gets away from $1,$ the
estimation of $\mu$ gets better and better while the latter indicates that the
estimation of $\sigma$ is not good when $\alpha$ is around $1.5$ but for
smaller values, it might be considered as acceptable. It is to be noted that,
in regards to the estimation of $\mu,$ the results are extremely poor when the
stability index is very close to $1.$ This may be explained by the fact that,
in this work, we only consider $\alpha$-values lying between $1$ and $2$ and
in this case the location parameter is equal to the distribution mean and is
estimated as such. When $\alpha$ is less than or equal to $1,$ the mean does
not exist and therefore the EVT-based estimation of $\mu$ is very bad when
$\alpha$ is near $1.$%

%

\begin{table}[h]  \centering
\begin{tabular}
[c]{cccccccc}%
$\alpha$ & \multicolumn{1}{|c}{$\sigma$} & \multicolumn{1}{|c}{$\widehat
{\alpha}_{n}$} & abs bias & mse & conf int & length & cov prob\\ \hline \hline
\multicolumn{1}{c|}{} & \multicolumn{1}{|l}{$0.1$} &
\multicolumn{1}{|c}{$1.114$} & $0.014$ & $0.011$ & $0.978,1.301$ & $0.323$ &
$0.91$\\
\multicolumn{1}{l}{$1.1$} & \multicolumn{1}{|l}{$0.5$} &
\multicolumn{1}{|c}{$1.110$} & $0.010$ & $0.011$ & $0.975,1.294$ & $0.319$ &
$0.91$\\
\multicolumn{1}{c|}{} & \multicolumn{1}{|l}{$1.0$} &
\multicolumn{1}{|c}{$1.114$} & $0.014$ & $0.012$ & $0.979,1.298$ & $0.319$ &
$0.90$\\ \hline
\multicolumn{1}{c|}{} & \multicolumn{1}{|l}{$0.1$} &
\multicolumn{1}{|c}{$1.227$} & $0.027$ & $0.013$ & $1.075,1.435$ & $0.359$ &
$0.91$\\
\multicolumn{1}{l}{$1.2$} & \multicolumn{1}{|l}{$0.5$} &
\multicolumn{1}{|c}{$1.217$} & $0.017$ & $0.013$ & $1.066,1.422$ & $0.356$ &
$0.91$\\
\multicolumn{1}{c|}{} & \multicolumn{1}{|l}{$1.0$} &
\multicolumn{1}{|c}{$1.231$} & $0.031$ & $0.014$ & $1.081,1.436$ & $0.354$ &
$0.89$\\ \hline
\multicolumn{1}{l}{} & \multicolumn{1}{|l}{$0.1$} &
\multicolumn{1}{|c}{$1.633$} & $0.133$ & $0.045$ & $1.421,1.931$ & $0.510$ &
$0.63$\\
\multicolumn{1}{l}{$1.5$} & \multicolumn{1}{|l}{$0.5$} &
\multicolumn{1}{|c}{$1.630$} & $0.130$ & $0.044$ & $1.419,1.923$ & $0.504$ &
$0.65$\\
\multicolumn{1}{l}{} & \multicolumn{1}{|l}{$1.0$} &
\multicolumn{1}{|c}{$1.627$} & $0.127$ & $0.046$ & $1.414,1.927$ & $0.513$ &
$0.65$\\ \hline
\multicolumn{1}{l}{} & \multicolumn{1}{|l}{$0.1$} &
\multicolumn{1}{|c}{$2.403$} & $0.603$ & $0.492$ & $2.013,3.004$ & $0.991$ &
$0.36$\\
\multicolumn{1}{l}{$1.8$} & \multicolumn{1}{|l}{$0.5$} &
\multicolumn{1}{|c}{$2.418$} & $0.618$ & $0.521$ & $2.032,3.009$ & $0.978$ &
$0.34$\\
\multicolumn{1}{l}{} & \multicolumn{1}{|l}{$1.0$} &
\multicolumn{1}{|c}{$2.405$} & $0.605$ & $0.509$ & $2.020,2.994$ & $0.974$ &
$0.36$\\ \hline
\multicolumn{1}{l}{} & \multicolumn{1}{l}{} &  &  &  &  &  &
\end{tabular}
\caption{Simulation results of the estimation of the shape parameter of a symmetric Lévy-stable distribution based on 1000 samples of 3000 observations.}\label{Tab1}%
\end{table}%
%

\begin{table}[h]  \centering
\begin{tabular}
[c]{cccccccc}%
$\alpha$ & \multicolumn{1}{|c}{$\sigma$} & \multicolumn{1}{|c}{$\widehat{\mu
}_{n}$} & abs bias & mse & conf int & length & cov prob\\ \hline \hline
\multicolumn{1}{l}{} & \multicolumn{1}{|l}{$0.1$} &
\multicolumn{1}{|c}{$0.025$} & $0.025$ & $0.035$ & $-0.432,0.482$ & $0.915$ &
$0.89$\\
\multicolumn{1}{l}{$1.2$} & \multicolumn{1}{|l}{$0.5$} &
\multicolumn{1}{|c}{$-0.058$} & $0.058$ & $0.213$ & $-2.158,2.042$ & $4.200$ &
$0.88$\\
\multicolumn{1}{l}{} & \multicolumn{1}{|l}{$1.0$} &
\multicolumn{1}{|c}{$0.171$} & $0.171$ & $0.542$ & $-2.878,3.220$ & $6.098$ &
$0.85$\\ \hline
\multicolumn{1}{l}{} & \multicolumn{1}{|l}{$0.1$} &
\multicolumn{1}{|c}{$0.001$} & $0.001$ & $0.002$ & $-0.094,0.096$ & $0.190$ &
$0.92$\\
\multicolumn{1}{l}{$1.3$} & \multicolumn{1}{|l}{$0.5$} &
\multicolumn{1}{|c}{$0.004$} & $0.004$ & $0.04$ & $-0.471,0.479$ & $0.950$ &
$0.92$\\
\multicolumn{1}{l}{} & \multicolumn{1}{|l}{$1.0$} &
\multicolumn{1}{|c}{$0.008$} & $0.008$ & $0.16$ & $-0.942,0.958$ & $1.900$ &
$0.92$\\ \hline
\multicolumn{1}{l}{} & \multicolumn{1}{|l}{$0.1$} &
\multicolumn{1}{|c}{$0.003$} & $0.003$ & $0.000$ & $-0.038,0.044$ & $0.083$ &
$0.94$\\
\multicolumn{1}{l}{$1.5$} & \multicolumn{1}{|l}{$0.5$} &
\multicolumn{1}{|c}{$-0.010$} & $0.010$ & $0.005$ & $-0.207,0.186$ & $0.393$ &
$0.94$\\
\multicolumn{1}{l}{} & \multicolumn{1}{|l}{$1.0$} &
\multicolumn{1}{|c}{$0.013$} & $0.013$ & $0.028$ & $-0.394,0.421$ & $0.814$ &
$0.94$\\ \hline
\multicolumn{1}{l}{} & \multicolumn{1}{l}{} &  &  &  &  &  &
\end{tabular}
\caption{Simulation results of the estimation of the location parameter of a symmetric Lévy-stable distribution based on 1000 samples of 3000 observations.}\label{Tab2}%
\end{table}%
%

\begin{table}[h] \centering
\begin{tabular}
[c]{cccccccc}%
$\alpha$ & \multicolumn{1}{|c}{$\sigma$} & \multicolumn{1}{|c}{$\widehat
{\sigma}_{n}$} & abs bias & mse & conf int & length & cov prob\\ \hline \hline
\multicolumn{1}{c|}{} & \multicolumn{1}{|l}{$0.1$} &
\multicolumn{1}{|c}{$0.113$} & $0.013$ & $0.002$ & $0.080,0.164$ & $0.150$ &
$0.81$\\
\multicolumn{1}{l}{$1.2$} & \multicolumn{1}{|l}{$0.5$} &
\multicolumn{1}{|c}{$0.570$} & $0.070$ & $0.047$ & $0.400,0.836$ & $0.750$ &
$0.80$\\
\multicolumn{1}{c|}{} & \multicolumn{1}{|l}{$1.0$} &
\multicolumn{1}{|c}{$1.128$} & $0.128$ & $0.164$ & $0.790,1.653$ & $1.506$ &
$0.81$\\ \hline
\multicolumn{1}{c|}{} & \multicolumn{1}{|l}{$0.1$} &
\multicolumn{1}{|c}{$0.119$} & $0.019$ & $0.002$ & $0.087,0.168$ & $0.081$ &
$0.70$\\
\multicolumn{1}{l}{$1.3$} & \multicolumn{1}{|l}{$0.5$} &
\multicolumn{1}{|c}{$0.609$} & $0.109$ & $0.138$ & $0.444,0.852$ & $0.408$ &
$0.70$\\
\multicolumn{1}{c|}{} & \multicolumn{1}{|l}{$1.0$} &
\multicolumn{1}{|c}{$1.190$} & $0.190$ & $0.218$ & $0.863,1.678$ & $0.815$ &
$0.69$\\ \hline
\multicolumn{1}{l}{} & \multicolumn{1}{|l}{$0.1$} &
\multicolumn{1}{|c}{$0.158$} & $0.058$ & $0.021$ & $0.100,0.255$ & $0.155$ &
$0.60$\\
\multicolumn{1}{l}{$1.5$} & \multicolumn{1}{|l}{$0.5$} &
\multicolumn{1}{|c}{$0.790$} & $0.290$ & $0.522$ & $0.500,1.276$ & $0.776$ &
$0.60$\\
\multicolumn{1}{l}{} & \multicolumn{1}{|l}{$1.0$} &
\multicolumn{1}{|c}{$1.664$} & $0.664$ & $1.442$ & $1.065,2.664$ & $1.599$ &
$0.60$\\ \hline
\multicolumn{1}{l}{} & \multicolumn{1}{l}{} &  &  &  &  &  &
\end{tabular}
\caption{Simulation results of the estimation of the scale parameter of a symmetric Lévy-stable distribution based on 1000 samples of 3000 observations.}\label{Tab3}%
\end{table}%

\end{document}